\documentclass[12pt]{article}
\usepackage{amsmath}
\usepackage{amssymb}
\usepackage{amscd}
\usepackage{amsthm}
\theoremstyle{plain}
\newtheorem{Thm}{Theorem}[section]
\newtheorem{Prob}[Thm]{Problem}
\newtheorem{Prop}[Thm]{Proposition}
\newtheorem{Conj}[Thm]{Conjecture}
\theoremstyle{definition}
\newtheorem{Defn}[Thm]{Definition}
\theoremstyle{remark}
\newtheorem{Rem}[Thm]{Remark}
\numberwithin{equation}{section}
\def\Exc{\mathop{\mathrm{Exc}}}
\def\Rat{\mathop{\mathrm{Rat}}}
\def\Bs{\mathop{\mathrm{Bs}}}
\def\Supp{\mathop{\mathrm{Supp}}}
\def\Hilb{\mathop{\mathrm{Hilb}}}
\title{On the projective fourfolds with almost numerically positive canonical divisors}
\author{ {Shigetaka FUKUDA}}
\date{\empty}

\begin{document}
\maketitle \thispagestyle{empty}
\pagestyle{myheadings}
\markboth{Shigetaka FUKUDA}{almost numerically positive}
\begin{abstract}
Let $X$ be a four-dimensional projective variety defined over the field of complex numbers with only terminal singularities.
We prove that if the intersection number of the canonical divisor $K$ with every very general curve is positive ($K$ is almost numerically positive) then every very general proper subvariety of $X$ is of general type in the viewpoint of geometric Kodaira dimension.
We note that the converse does not hold for simple abelian varieties.
\end{abstract}

\section{Introduction}

In this paper, every algebraic variety is defined over the field of complex numbers.
We follow the standard notation and terminology of the minimal model theory.
Let $X$ be a projective variety.

\begin{Defn}
The geometric Kodaira dimension $\kappa_{geom} (X)$ denotes the Kodaira dimension of non-singular models of $X$.
A $\mathbf{Q}$-Cartier divisor $D$ on $X$ is {\it almost numerically positive} ({\it almost nup}, for short), if there exists a union $F$ of at most countably many prime divisors on $X$ such that $(D , C) > 0$ for every curve $C \nsubseteq F$ (i.e.\ if $(D , C) > 0$ for every very general curve $C$).
We say that $D$ is {\it quasi-numerically positive} ({\it quasi-nup}, for short), if $D$ is nef and almost nup.
\end{Defn}

\begin{Defn}
For a set $P$, a subset $S$ of $P$ is {\it covered} by subsets $S_i$ of $P$, if $S \subset \bigcup_i S_i$.
\end{Defn}

Of course, the quasi-nupness implies the almost nupness.
We note that the abbreviation ``nup'' for the term ``numerically positive'', which means being strictly positive in intersection number with every curve, has been used by specialists.

Recently, Ambro (\cite{Am}) reduced the famous abundance conjecture that the nef canonical divisor $K_X$ on every minimal variety $X$ should be semi-ample, in dimension four, to the following:

\begin{Prob}
Let $X$ be four-dimensional and with only terminal singularities.
When the canonical divisor $K_X$ is quasi-nup, is $K_X$ semi-ample?
\end{Prob}

We look at the problem from the birational-geometric standpoint and state the following conjecture.

\begin{Conj}
Let $X$ be $n$-dimensional and with only terminal singularities.
The canonical divisor $K_X$ is almost nup if and only if $X$ is of general type $($i.e.\ $\kappa (X) = n$$)$.
\end{Conj}

This holds in dimension $n \leq 3$ by virtue of the minimal model and the abundance theorems (see Proposition \ref{Prop:Ideal}).
Motivated by the conjecture, we investigate the varieties with almost nup canonical divisors and obtain the following:

\begin{Thm}[Main Theorem]\label{Thm:MT}
Let $X$ be four-dimensional and with only terminal singularities.
If $K_X$ is almost nup $($almost numerically positive$)$, then the locus $\bigcup \{ D ; \thickspace D$ is a closed subvariety $\subsetneqq X$ and $\kappa_{geom} (D) < \dim D \} $ can be covered by at most countably many prime divisors on $X$.
\end{Thm}

This means that the fourfolds with almost nup canonical divisors have the geometric property that every very general proper subvariety is of general type (a kind of weak hyperbolicity).
The proof, which is given in Section \ref{Sect:PMT}, depends on the facts: some stability of the almost nupness under pulling-back, the countability of the components of the Hilbert scheme and the deformation invariance of the Kodaira dimension due to Tsuji (\cite{Ts}) and Siu (\cite{Si}).

The converse statement of Main Theorem \ref{Thm:MT} does not hold, because all proper subvarieties of simple abelian varieties are of general type (Ueno \cite{Ue}, Corollary 10.10).
So we give the following formulation, inspired by Lang's conjectures (\cite{La}, Section IV.5) concerning the concept of measure hyperbolicity.

\begin{Thm}\label{Thm:CMT}
Let $X$ be four-dimensional and with only terminal singularities.
The canonical divisor $K_X$ is almost nup if and only if the locus $\bigcup \{ D ; \thickspace D$ is a closed subvariety of $X$ ($D \subsetneqq X$ or $D=X$), $D$ is a rational curve or is birationally equivalent to a minimal variety with numerically trivial canonical divisor$\}$ can be covered by at most countably many prime divisors on $X$.
\end{Thm}

The proof is given in Section \ref{Sect:PCMT}, by using Shokurov's existence of fourfold flips (\cite{Sh}).

\begin{Rem}
The following formulation for surfaces is derived from Lang (\cite{La}, Section IV.5):

Let $S$ be a non-singular projective surface.
The canonical divisor $K_S$ is almost nup if and only if the locus $\bigcup \{ D ; \thickspace D$ is a closed subvariety of $S$ ($D \subsetneqq S$ or $D=S$), $D$ is a rational curve or is birationally equivalent to a simple abelian variety$\}$ can be covered by at most countably many prime divisors on $S$.

This is due to the fact that all algebraic surfaces with numerically trivial canonical divisor except simple abelian surfaces can be covered by a family of elliptic curves (from classical results for abelian, bielliptic and Enriques surfaces and from the result (\cite{MoMu}, the appendix) of Bogomolov-Mori-Mukai-Mumford for K3 surfaces). 
\end{Rem}

\section{The Proof of Main Theorem \ref{Thm:MT}}\label{Sect:PMT}

\begin{Prop}\label{Prop:Ideal}
Let $X$ be $n$-dimensional and with only terminal singularities.
Assume that the minimal model and the abundance conjectures hold in dimension $n$.
Then the following three conditions are equivalent:

$(1)$ $X$ is of general type $($i.e.\ $\kappa (X) = n$$)$.

$(2)$ There exists an effective reduced divisor $F$ such that $(K_X , C) > 0$ for every curve $C \nsubseteq F$.

$(3)$ $K_X$ is almost nup.
\end{Prop}

\begin{proof}
We divide the situation into three cases, by the value of the Kodaira dimension $\kappa (X)$. 

{\it The case where} $\kappa (X) = n$.
There exist a non-singular projective variety $Y$ and a birational morphism $f:Y \to X$ such that $\vert f^* (mK_X) \vert = \vert D \vert + E$ for some positive integer $m$, where $m K_X$ is Cartier, $\Bs \vert D \vert = \emptyset$, the morphism $\Phi_{\vert D \vert}$ is birational and $E \geq 0$.
Hence $(f^* (mK_X) , C') > 0$ for every curve $C' \nsubseteq \Supp (E) \cup \Exc (\Phi_{\vert D \vert})$.
Consequently $(mK_X , C) > 0$ for every curve $C \nsubseteq f(\Supp (E) \cup \Exc (\Phi_{\vert D \vert}) \cup \Exc (f) )$.
This implies (2) and (3). 

{\it The case where} $0 \leq \kappa (X) \leq n-1$. 
There exist a non-singular projective variety $Y$, a projective variety $Z$ with only terminal singularities and birational morphisms $f:Y \to X$ and $g:Y \to Z$ such that $K_Z$ is nef, where $m K_Z$ is Cartier, $\Bs \vert m K_Z \vert = \emptyset$ and the morphism $\Phi_{\vert m K_Z \vert}$ is with only connected fibers for some $m \geq 1$.
Let $F$ be a general fiber of $\Phi_{\vert m K_Z \vert} \cdot g$ (note that $1 \leq \dim F \leq n$). 
We have $(K_Y - g^* K_Z , C' ) = 0$ for $C' = \bigcap_{i=1}^{\dim F -1} g^{-1} (H_i) \cap F$ where $H_i$ is a general hyperplane section of $Z$, because $(K_Y - g^* K_Z) \vert_F$ is $(g \vert_F)$-exceptional.
Hence $(K_Y , C' ) = 0$.
Therefore, from the fact that $K_Y \geq f^* K_X $ and from the choice of $C'$, we have $(K_X, f(C')) \leq 0$.
This means that the conditions (2) and (3) do not hold.

{\it The case where} $\kappa (X) = -\infty$.
In this case $X$ is uniruled.
Thus, by Miyaoka-Mori (\cite{MiMo}), for every general point $x \in X$, there exists a curve $C$ such that $x \in C$ and $(K_X , C) < 0$.
This means that the conditions (2) and (3) do not hold.
\end{proof}

\begin{Prop}\label{Prop:Real}
Let $X$ be $n$-dimensional and with only terminal singularities.
Assume that the minimal model and the abundance conjectures hold in dimension $\leq n-1$.
If $K_X$ is almost nup, then the locus $\bigcup \{ D ; \thickspace D$ is a closed subvariety $\subsetneqq X$ and $\kappa_{geom} (D) < \dim D \}$ can be covered by at most countably many prime divisors on $X$.
\end{Prop}

\begin{proof}
Assuming that $K_X$ is almost nup and that however the locus $\bigcup \{ D ; \thickspace D$ is a closed subvariety $\subsetneqq X$ and $\kappa_{geom} (D) < \dim D \}$ cannot be covered by at most countably many prime divisors on $X$, we derive a contradiction.

Let $\mathcal{H} \subset X \times \Hilb (X)$ be the universal family parametrized by the Hilbert scheme $\Hilb (X)$.
By the countability of the components of $\Hilb (X)$, we have an irreducible component $V$ of $\mathcal{H}$ with surjective projection morphisms $f:V \to X$ and $g:V \to T ( \subset \Hilb (X))$ from $V$ to projective varieties $X$ and $T$ respectively, such that $f(g^{-1} (t)) \subsetneqq X$ for every $t \in T$ and that the locus $\bigcup \{ D ; \thickspace D ( \subsetneqq X)$ is a closed subvariety, $\kappa_{geom} (D) < \dim D$ and $D = f(g^{-1} (t))$ for some $t \in T \}$ cannot be covered by at most countably many prime divisors on $X$.
Let $\nu : V_{norm} \to V$ be the normalization.
We consider the Stein factorization of $g \nu$ into the finite morphism $g_1 : S \to T$ from a projective normal variety $S$ and the morphism $g_2 : V_{norm} \to S$ with an algebraically closed extension $\Rat V / \Rat S$.
Put $V^* :=$ [the image of the morphism $(f \nu, g_2) : V_{norm} \to X \times S$].
\begin{equation}
 \begin{CD}
V_{norm} @>>> V^* @>\text{embedding}>> X \times S @>>> S    \\
@.   @VVV @VVV @VV\text{$g_1$}V \\
@.   V   @>\text{embedding}>> X \times T @>>> T    \\
@.   @.   @VVV \\
@.   @.   X
 \end{CD}
\end{equation}
Note that every fiber of the morphism $g:V \to T$ consists of a finite number of fibers of the projection morphism from $V^*$ to $S$.
Thus we may replace $(V,T)$ by $(V^* , S)$ and assume that the extension $\Rat V / \Rat T$ is algebraically closed.

From the deformation invariance of the Kodaira dimension due to Tsuji (\cite{Ts}) and Siu (\cite{Si}), $\kappa_{geom} (g^{-1} (t))$ for general $t \in T$ is constant.

In the case where $\kappa_{geom} (g^{-1} (t)) = \dim g^{-1} (t)$ for general $t \in T$, there exists a subvariety $T_0 \subsetneqq T$ such that the locus $\bigcup \{ D ; \thickspace D ( \subsetneqq X)$ is a closed subvariety, $\kappa_{geom} (D) < \dim D$ and $D = f(g^{-1} (t))$ for some $t \in T_0 \}$ cannot be covered by at most countably many prime divisors on $X$.
Thus we can replace $(V, T)$ by $(V_1 , T_1)$, where $V_1$ and $T_1$ are projective varieties such that $V_1$ is some suitable irreducible component of $g^{-1} (T_0)$ and $T_1 = g(V_1)$ .
Because $\dim V_1 < \dim V$, by repeating this process of replacement, we can reduce the assertion to the next case.

Now consider the case where $\kappa_{geom} (g^{-1} (t)) < \dim g^{-1} (t)$ for general $t \in T$.
If $\dim V > \dim X$, then $\dim f^{-1} (x) \geq 1$ for all $x \in X$, thus $\dim g(f^{-1} (x)) \geq 1$ (this means that $g(f^{-1} (x)) \cap H \ne \emptyset$ for every hyperplane section $H$ of $T$).
Therefore, by repeating the process of cutting $T$ by general hyperplanes, we can reduce the assertion to the subcase where $\dim V = \dim X$.
Take birational morphisms $\alpha : X' \to X$ and $\beta : V' \to V$ from non-singular projective varieties with a generically finite morphism $\gamma : V' \to X'$ such that $\alpha \gamma = f \beta$.
\begin{equation}
 \begin{CD}
V @<\text{$\beta$}<< V' \\
@V\text{$f$}VV @VV\text{$\gamma$}V \\
X @<<\text{$\alpha$}< X'
 \end{CD}
\end{equation}
Because $K_X$ is almost nup and $K_{X'} - \alpha^* K_X$ is effective, $K_{X'}$ is almost nup.
Here we have $K_{V'} = \gamma^{*} K_{X'} + R_{\gamma}$, where $R_{\gamma}$ is the ramification divisor (which is effective) for $\gamma$.
Thus $K_{V'}$ becomes almost nup.
For a very general fiber $F' = (g \beta )^{-1} (t) $ for $g \beta$ (i.e.\ the point $t$ does not belong to some fixed union of at most countably many prime divisors on $T$), also $K_{F'}$ is almost nup.
Thus from the assumption of the proposition and from Proposition \ref{Prop:Ideal}, $\kappa (F') = \dim F'$.
So $\kappa_{geom} (F) = \dim F$ for every general fiber $F$ of $g$.
This is a contradiction!
\end{proof} 

At last we prove the main theorem.

\begin{proof}[Proof of Main Theorem \ref{Thm:MT}]
Proposition \ref{Prop:Real} implies the assertion, because the minimal model and the abundance conjectures hold in dimension $\leq 3$ (Kawamata \cite{Ka}, Miyaoka \cite{Mi}, Mori \cite{Mo}).
\end{proof}

\section{The Proof of Theorem \ref{Thm:CMT}}\label{Sect:PCMT}

\begin{Prop}\label{Prop:Prep}
Let $X$ be $n$-dimensional and with only terminal singularities.
If $K_X$ is almost nup, then $X$ is not birationally equivalent to any minimal variety with numerically trivial canonical divisor.
\end{Prop}

\begin{proof}
Assuming that $K_X$ is almost nup and that $X$ is birationally equivalent to a minimal variety $M$ with numerically trivial canonical divisor $K_M$, we derive a contradiction.
We take a common resolution $\mu : Y \to X$ and $\nu : Y \to M$.
Then $K_Y (\geq \mu^*K_X)$ is almost nup.

Let $H$ be a hyperplane section of $M$.
We have that $(K_Y, (\nu^* H)^{n-1} ) = (K_M, H^{n-1}) = 0$ because $K_Y - \nu^* K_M$ is $\nu$-exceptional.
Thus $K_Y$ is not almost nup.
This is a contradiction!
\end{proof}

\begin{Prop}\label{Prop:Conv}
Let $X$ be $n$-dimensional and with only terminal singularities.
Assume that the minimal model conjecture holds in dimension $n$ and the log minimal model and the log abundance conjectures hold in dimension $\leq n-1$.
If the locus $\bigcup \{ D ; \thickspace D$ is a closed subvariety of $X$ ($D \subsetneqq X$ or $D=X$), $D$ is a rational curve or is birationally equivalent to a minimal variety with numerically trivial canonical divisor$\}$ can be covered by at most countably many prime divisors on $X$, then $K_X$ is almost nup.
\end{Prop}

\begin{proof}
We obtain a $\mathbf{Q}$-factorialization $\mu: X' \to X$ such that $K_{X'} = \mu^* K_X$, by running the relative minimal model program for some desingularization of $X$.
If $X'$ has no minimal model, then some nonempty Zariski open subset of $X'$ can be covered by rational curves.
Thus we may assume that $X'$ has a minimal model.
So let $X''$ be a minimal model that is an output from the minimal model program for $X'$.
Here we have a common resolution $\alpha : X_0 \to X'$ and $\beta : X_0 \to X''$ such that $\alpha^* K_{X'} = \beta^* K_{X''} + E$ where $E$ is an effective $\beta$-exceptional divisor.
\begin{equation}
 \begin{CD}
X @<\text{$\mu$}<< X' @<\text{$\alpha$}<< X_0 \\
@.   @.   @VV\text{$\beta$}V \\
@.   @.   X''
 \end{CD}
\end{equation}
Therefore, if $K_X$ is not almost nup, then $K_{X''}$ is not quasi-nup, hence $K_{X''}$ is not big but semi-ample from Ambro (\cite{Am}).

This means that, if $K_X$ is not almost nup and $X$ is not birationally equivalent to any minimal fourfold with numerically trivial canonical divisor, then some nonempty Zariski open subset of $X$ is covered by proper subvarieties with geometric Kodaira dimension $\leq 0$.
\end{proof}

\begin{proof}[Proof of Theorem \ref{Thm:CMT}]
Main Theorem \ref{Thm:MT} and Proposition \ref {Prop:Prep} imply the ``only if'' part.

Proposition \ref{Prop:Conv} implies the ``if'' part, because the minimal model conjecture holds in dimension four (Shokurov \cite{Sh} for the existence of flips and Kawamata-Matsuda-Matsuki (\cite{KMM}, Theorem 5.1.15) for the termination of flips) and the log minimal model and the log abundance conjectures hold in dimension $\leq 3$ (Kawamata \cite{Ka2} and Fujita \cite{Ft} for surfaces, Shokurov \cite{Sh3} for three-dimensional log flips and Keel-Matsuki-McKernan \cite{KeMaMc} for the three-dimensional log abundance).
\end{proof}

\bigskip
Faculty of Education, Gifu Shotoku Gakuen University

Yanaizu-cho, Gifu 501-6194, Japan

fukuda@ha.shotoku.ac.jp

\end{document}